\documentclass[11pt]{article}

\usepackage[top=1.5in, bottom=1.5in, left=1in, right=1in]{geometry}
\usepackage{microtype} 

\usepackage{amsmath}   
\usepackage{amsfonts}
\usepackage{amssymb}
\usepackage{amsthm}
\usepackage{mathrsfs}
\usepackage{faktor}
\usepackage{hyphenat}
\usepackage{graphicx}
\usepackage[color,matrix,arrow,all]{xy} 
\usepackage{tikz}
\usepackage{tikz-cd}
\usetikzlibrary{matrix,calc}
\usepackage{eepic} 

\usepackage{appendix}
\usepackage{color}

\newcommand{\dom}{\mathrm{Dom}}

\newcommand{\la}{\langle}
\newcommand{\ra}{\rangle}

\newcommand{\bi}{\begin{itemize}}
\newcommand{\ei}{\end{itemize}}
\newcommand{\beq}{\begin{equation}}
\newcommand{\eeq}{\end{equation}}

\DeclareMathOperator{\GM}{GM}
\DeclareMathOperator{\RLM}{RLM}
\DeclareMathOperator{\RM}{RM}
\DeclareMathOperator{\CM}{CM}

\theoremstyle{plain}
\newtheorem{T}{Theorem}[section]
\newcommand{\bt}{\begin{T}}
\newcommand{\et}{\end{T}}

\newtheorem{Proposition}[T]{Proposition}
\newcommand{\bp}{\begin{Proposition}}
\newcommand{\ep}{\end{Proposition}}

\newtheorem{Lemma}[T]{Lemma}
\newcommand{\bl}{\begin{Lemma}}
\newcommand{\el}{\end{Lemma}}

\newtheorem{Corol}[T]{Corollary}
\newcommand{\bc}{\begin{Corol}}
\newcommand{\ec}{\end{Corol}}

\newtheorem{Result}[T]{Result}
\newcommand{\br}{\begin{Result}}
\newcommand{\er}{\end{Result}}

\theoremstyle{definition}
\newtheorem{Example}[T]{Example}
\newcommand{\be}{\begin{Example}}
\newcommand{\ee}{\end{Example}}

\newtheorem{Problem}[T]{Problem}
\newcommand{\bq}{\begin{Problem}}
\newcommand{\eq}{\end{Problem}}

\newtheorem{Remark}[T]{Remark}
\newcommand{\brm}{\begin{Remark}}
\newcommand{\erm}{\end{Remark}}

\newtheorem{Discussion}[T]{Discussion}
\newcommand{\bds}{\begin{Discussion}}
\newcommand{\eds}{\end{Discussion}}

\newtheorem{Definition}[T]{Definition}
\newcommand{\bd}{\begin{Definition}}
\newcommand{\ed}{\end{Definition}}

\newtheorem{Construction}[T]{Construction}
\newcommand{\bco}{\begin{Construction}}
\newcommand{\eco}{\end{Construction}}

\title{Aperiodic Flows on Finite Semigroups II: Further Developments}
\author{Stuart Margolis \and John Rhodes}

\title{Aperiodic Flows on Finite Semigroups II: Smallish Monoids Suffice for Complexity 1 }
\author{Stuart Margolis and John Rhodes}

\date{\today
\footnote{MSC classes: 20M10, 20M20, 20M30, 20M35}
}

\begin{document}

\maketitle

\begin{abstract}

A smallish monoid $M$ is a monoid that has a unique 0-minimal ideal $I(M)$ that is a 0-simple subsemigroup and such that its regular $\mathcal{J}$-classes are the group of units and the two in $I(M)$. We show constructively how to embed an arbitrary finite semigroup $S$ into the evaluation semigroup of a smallish monoid $S^{Ev}$. We use the theory of flows to show that a group  mapping semigroup $S$ admits an aperiodic flow if and only if $S^{Ev}$ admits one. This reduces the computation of Krohn-Rhodes complexity 1 to the class of smallish monoids.
    
\end{abstract}
\section{Introduction}

This paper is a continuation of the paper \cite{FlowsI}. We will use all the notation and results from that paper. We summarize some of that information and other background material in the Appendices of this paper. 

The main result of this paper shows how to embed any semigroup into the evaluation semigroup (See Section 4 of \cite{Trans} or Appendix \ref{Eval}) of a $\GM$ semigroup with a very restricted structure.  More precisely, a {\em smallish monoid} $S$ is a monoid that has a unique 0-minimal ideal $I(S)$ and whose regular $\mathcal{J}$-classes are its group of units $U(S)$, $I(S) - \{0\}$ and $\{0\}$. Regular monoids with this property are called small monoids and have appeared in many places in semigroup theory and its applications.

The proof is constructive. To every $\GM$ semigroup $S$ we construct a smallish $\GM$ monoid $S^{Ev}$ such that $S$ embeds into the evaluation semigroup $Eval(S^{Ev})$. Furthermore, we use the theory of flows to show that $S$ has an aperiodic flow if and only if $S^{Ev}$ has an aperiodic flow. In particular, if it is known by induction that $\RLM(S)c \leq 1$, then $Sc = 1$ if and only if $S^{Ev}c=1$. This reduces computing complexity 1 to smallish monoids.

Although complexity has been proved decidable for all finite semigroups \cite{complexity1, complexityn}, the construction of $S^{Ev}$ is of independent interest. Furthermore, it has been known for more than 50 years that complexity of small monoids is decidable \cite{2J}. Thus our main result shows that adding non-regular $\mathcal{J}$-classes between the group of units and the unique 0-minimal ideal to small monoids defines a class general enough to compute complexity 1 for all semigroups.

\section{Preliminaries}

All semigroups are finite in this paper. Complexity means Krohn-Rhodes complexity. The complexity $Sc$ of a finite semigroup $S$ is the least integer $n$ such that $S$ divides a wreath product $$S \prec A_{0} \wr G_{1}\ldots A_{n-1}\wr G_{n}\wr A_{n}$$
where all the $A_{i}$ are aperiodic semigroups and all the $G_{j}$ are groups. The Krohn-Rhodes Theorem \cite{qtheory} guarantees that every finite semigroup $S$ has a well-defined complexity $Sc$. For background in complexity, see \cite[Chapter 4]{qtheory} and \cite{Arbib}.  
For all undefined notation see \cite{qtheory}, \cite{Arbib}, \cite{Eilenberg}. For background on the theory of flows, see \cite{Trans, FlowsI} or Appendix \ref{Flows}.

\medskip
\noindent
Examples are given in the form of a Rees matrix semigroup $\mathcal{M}^{0}(G,A,B,C)$ and additional generators. For historic reasons, instead of $C$ as above we use the transpose
$C^{T}: A \times B \to G^0$. Generators are given (usually not a minimal set) of the form 
$(a,g,b) \in \mathcal{M}^{0}(G,A,B,C)$ that generate all of $\mathcal{M}^{0}(G,A,B,C)$, and a
collection of elements of the monoid $\RM_{B}(G)$ of $|B| \times |B|$ row-monomial matrices with entries in $G^{0}$. We write an element $X$ of $\RM_{B}(G)$ as a $G$-labeled partial function. The edge $i \rightarrow gj, i,j \in B, g \in G$ means that the $(i,j)$ entry of $X$ is $g$. We identify $X$ with the set of such edges, one for each non-zero row of $X$. 

We have the dual representation for the monoid of $|A| \times |A|$ column-monomial matrices $\CM_{A}(G)$ with entries in $G$. We assume that each generator $X \in \RM_{B}(G)$ is linked with at least one $Y \in \CM_{A}(G)$. This means that $XC=CY$. The collection of all such pairs $(X,Y)$ is called the translational hull of $\mathcal{M}^{0}(G,A,B,C)$ \cite[Chapter 5]{qtheory}. We do not need any facts about translational hulls except that it guarantees that our semigroups have $\mathcal{M}^{0}(G,A,B,C)$ as their unique 0-minimal ideals. We do not give a $Y$ such that $(X,Y)$ is in the translational hull of $\mathcal{M}^{0}(G,A,B,C)$. The reader can compute that directly. 

A generalized group-mapping (written $\operatorname{GGM}$) semigroup $S$ is a semigroup that has a unique 0-minimal ideal $I(S)$ that is also a 0-simple subsemigroup $I(S)=\mathcal{M}^{0}(G,A,B,C)$ such that $S$ acts faithfully on both the left and the right of $I(S)$ by multiplication. $S$ is a group-mapping semigroup (written $\GM$) if it is a $\operatorname{GGM}$ semigroup such that the maximal subgroup of $I(S)$ is non-trivial. An important property of $\GM$ semigroups is that every finite semigroup $S$  has a $\GM$ image that has the same complexity as $S$.


\medskip
\noindent
With the notation above, we have a transformation semigroup $(G \times B,S)$. If $s \in S$ is given as above by a $G$-labeled partial function, 
we have $(g,b)s=(gh,b')$ if $b \rightarrow hb'$ is an edge of $s$ and undefined if there is no edge beginning with $b$ in $s$.

\section{Embedding arbitrary semigroups into Evaluation Semigroups}

Proposition 5.16 of \cite{Trans} proves that any $\GM$ semigroup $S$ embeds into its evaluation semigroup $Eval(S)$. See \cite{Trans} or Appendix \ref{Eval} of this paper for background on $Eval(S)$. Since $Z_{2} \times PT_{n}$ is a $\GM$ semigroup, where $PT_n$ is the monoid of all partial functions on an $n$-set, it follows that every semigroup embeds into the evaluation semigroup of some $\GM$ semigroup. The purpose of this section is to improve this result. 

Recall that a small monoid is a monoid $M$ of the form $U(M) \cup I(M)$ where $U(M)$ is the group of units of $M$ and $I(M) =\mathcal{M}^{0}(G,A,B,C)$ is its unique 0-minimal ideal that is also a 0-simple semigroup. In other words, small monoids are precisely the regular monoids with exactly 2 non-zero $\mathcal{J}$-classes. We generalize this definition.

\bd

A monoid $M$ is called a {\em smallish monoid} if and only if the ideal $M - U(M)$ of non-units of $M$ is a nilpotent ideal extension of a 0-simple semigroup. That is, if $I$ is the ideal of non-units $M-U(M)$ of $M$, then there is a $k>0$ such that $I^{k}$ is the unique 0-minimal ideal of $M$ and is a 0-simple semigroup.

\ed

The following Lemma justifies the name smallish monoids for this class.

\bl \label{smallcond}

A monoid $M$ is a smallish monoid if and only if $M$ has a unique  0-minimal ideal $I(M)$ that is a regular subsemigroup of $M$ and has exactly 2 non-zero regular $\mathcal{J}$-classes, namely $I(M) -\{0\}$ and its group of units $U(M)$.  
\el

\begin{proof}

Let $M$ be a smallish monoid and let $I=M -U(M)$. Then there is a $k>0$ such that $I^{k}=I(M)$ where $I(M)$ is an ideal that is also a 0-simple semigroup. It follows that every idempotent of $I$ is in $I(M)$ and this implies that the only regular $\mathcal{J}$-classes are $U(M), I(M)-\{0\}$ and $0$. Therefore, $M$ has 2 non-zero regular $\mathcal{J}$-classes and  $I(M)$ is the unique 0-minimal ideal of $M$.

Conversely assume that $M$ is a monoid with a unique 0-minimal ideal $I(M)$ that is a 0-simple subsemigroup of $M$ and exactly 2 non-zero regular $\mathcal{J}$-classes. Clearly the regular non-zero $\mathcal{J}$-classes are the group of units $U(M)$ and $I(M)-\{0\}$. Let $I=M -U(M)$. It follows that all idempotents in $I$ are in $I(M)$. Since it is well known that the unique idempotent power $I^{k}$ of $I$ satisfies $I^{k}=IE(I)I$ where $E(I)$ is the set of idempotents in $I$, it follows that $I$ is a nilpotent ideal extension of the 0-simple semigroup $I(M)$.

\end{proof}
If $M$ is a smallish monoid then $U(M) \cup I(M)$ is a small submonoid called the {\em small submonoid} 
$\operatorname{Sm}(M)$ of $M$. The following result follows from the Depth Decomposition Theorem. \cite{TilsonXII}.

\bp
Let $S$ be a smallish monoid. Then $Sc \leq 2$.
\ep

\be We give an example  of a smallish monoid $MTF$, called the Modified Tall Fork. This is a modification of Example 1 of \cite{MasterList}, the Tall Fork. See also Section 4.14 of \cite{qtheory}.


\noindent  $I(MTF)=\mathcal{M}^{0}(Z_{2},\{a_{1},\ldots , a_{7}\},\{1',3',1,2,3,4\},C)$  where \ $Z_{2} = \{\pm 1\}$ is the group of order 2 and the transpose of the structure matrix of $I(MTF)$ is given by: 

\bigskip

$C^{T} \ = \ \ $
\begin{tabular}{r|c c|r r r r|}
      &$1'$&$3'$& 1 & 2 & 3 & 4  \\ \hline
$a_1$ & 1 & 1   & 0 & 0 & 0 & 0  \\ 
$a_2$ & 1 & 0   & 0 & 0 & 0 & 0  \\ 
$a_3$ & 0 & 1   & 0 & 0 & 0 & 0  \\ \hline 
$a_4$ & 0 & 0   & 1 & 1 & 0 & 0  \\ 
$a_5$ & 0 & 0   & 0 & 1 & 1 & 0  \\ 
$a_6$ & 0 & 0 & 0 & 0 & 1 & 1  \\ 
$a_7$ & 0 & 0  & 1 & 0 & 0 & 1  \\ \hline
\end{tabular}

\bigskip

\noindent The generators of $MTF$ are all the elements of $I(MTF)$ together with the following elements of ${\rm RM}(\{1',3',1,2,3,4\},Z_{2})$. 

$\sigma = (1' \ 3')(1 \ 2 \ 3 \ 4)$ 


$r = \{(1' \to 1), \ (3' \to -3)\}$ 


\noindent The group of units $U(MTF)$ is the subgroup generated by $\sigma$ and is isomorphic to the cyclic group $Z_4$ of order 4. The poset of $\mathcal{J}$-classes is as follows. 

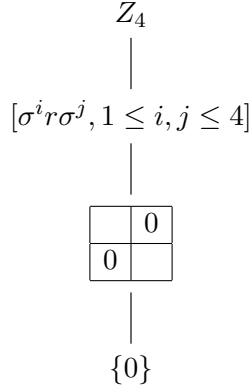
\begin{center}
\begin{tikzcd}[row sep=normal, ampersand replacement=\&]
Z_4 \arrow[d, dash] \\
{[\sigma^{i}r\sigma^{j}, 1 \leq i,j \leq 4]} \arrow[d, dash] \\
\begin{array}{|c|c|}
\hline
  & 0 \\ \hline
0 &   \\ \hline
\end{array} \arrow[d, dash] \\
\{0\}
\end{tikzcd}
\end{center}

\noindent where the matrix above $\{0\}$ represents the non-zero $\mathcal{J}$-class of $I(MTF)$ and has maximal subgroup $Z_{2}$. We note that the $\mathcal{J}$-class below $Z_{4}$ is a $4 \times 4$ null-$\mathcal{J}$-class with trivial Sch\"utzenberger group.

\noindent Therefore $MTF$ is a smallish monoid. By modifying the proof of Example 1. of \cite{MasterList} we have $MTFc=2$. On the other hand, it follows from Tilson's Theorem of \cite{2J} (see also Theorem \ref{2J} of this paper) that the small submonoid $Sm(MTF)$ has complexity 1. We will see that this example is generic as we see in the next theorem, that is the main theorem of this paper.

\ee
\bt\label{evalembed}

Let $(G \times B,S)$ be a $\GM$ transformation semigroup with 0-minimal ideal $\mathcal{M}^{0}(G,A,B,C)$ where $|B| = n$. Then we can construct another $\GM$ transformation semigroup $(G \times B^{Ev},S^{Ev})$ with the following properties.

\begin{enumerate}

\item{$S^{Ev}$ is a smallish $\GM$ monoid whose group of units $U(S^{Ev})$ is the cyclic group $Z_{n+1}$.}

\item{$S$ embeds into the evaluation semigroup $Eval(S^{Ev})$.}

\item{The small submonoid $Sm(S^{Ev})$ has complexity 1.}

\item{$S$ admits an aperiodic flow if and only if $S^{Ev}$ admits an aperiodic flow.}

\item{If $\RLM(S)c \leq 1$ then $Sc=1$ if and only if $S^{Ev}c = 1$.}

\end{enumerate}

\et 

\bc
The problem of deciding whether a semigroup $S$ has complexity 1 reduces to the case where $S$ is a smallish $\GM$ monoid.
\ec

\begin{proof}

It is known that for every semigroup $T$ there is a $\GM$ image $S$ such that $Sc=Tc$ \cite{Arbib, qtheory}. The result now follows from Theorem \ref{evalembed} Item 3. by induction on the cardinality of $S$. Here we use the fact that if $S$ is a $\GM$ semigroup then the cardinality of $\RLM(S)$ is strictly less than the cardinality of $S$.

\end{proof}

We prove Theorem \ref{evalembed} via the following construction of $(G \times B^{Ev},S^{Ev})$ from $(G \times B,S)$.

\bco

Let $S$ be a $\GM$ semigroup with $I(S) = \mathcal{M}^{0}(G,A,B,C)$. Let $B=\{1,\ldots , n\}$. We choose and fix an $a \in A$ and can assume that $C(1,a)=1$, so that $G$ is isomorphic to the $\mathcal{H}$-class $\{(a,g,1)|g \in G\}$ via the isomorphism that sends $g$ to $(a,g,1)$. We identify $G \times B$ with $\{a\}
 \times G \times B$. We can also assume that $C(i,a) \in \{0,1\}$ for $i=1, \ldots, n$. Such a normalization is always possible \cite[Chapter 7]{Arbib}, \cite[Chapter 4]{qtheory}.

Let $B^{Ev}=\{1,\ldots n\} \times \{0,\ldots n\}=B \times \{0,\ldots n\}$. 
We define $S^{Ev}$ to be the subsemigroup of 
$\RM_{B^{Ev}}(G)$ generated by the following elements of $\RM_{B^{Ev}}(G)$.

\begin{itemize}

\item{$t\colon (i,j) \mapsto (i,j+1)$ where we take $j+1$ modulo $n+1$.}

\item{Let $x \in S$ with $i \mapsto g_{i,x}\cdot ix$ for $i \in Dom(x), g_{i,x} \in G$. Define $h_{x}$ to have domain $\{(i,i)|i \in Dom(x)\}$ and defined by  $(i,i) \mapsto g_{i,x}\cdot (ix,0)$ for $(i,i) \in Dom(h_{x})$.}

\end{itemize}

\eco

We then have the transformation semigroup $(G \times B^{Ev},S^{Ev})$ defined by the action of $S^{Ev}$ on $G \times B^{Ev}$. We give the action of the generators of $S^{Ev}$ on $G \times B^{Ev}$.

\begin{itemize}

\item{($(g,(i,j))t=(g,(i,j+1)), g \in G, (i,j) \in B^{Ev}$.}

\item{$(g,(i,i))h_{x}=(gg_{i,x},(ix,0))$ for $x \in S, i \in Dom(x)$ with $ix =g_{i,x}\cdot ix$.}
    
\end{itemize}

We note the following important properties of $S^{Ev}$ that the reader can easily verify.

\begin{enumerate}

\item{The group of units $U(S^{Ev})$ of $S^{Ev}$ is generated by $t$ and is isomorphic to $Z_{n+1}$. For $i = 1, \ldots , n$ let $\mathcal{O}_{i}=\{(i,j)|j=0 \ldots , n\}$. Then the orbits of $U(S^{Ev})$ on $G \times B^{Ev}$ are the sets $ \{g\} \times \mathcal{O}_{i}, i = 1 , \ldots , n$ for $g \in G$.}


\item{Each element $h_{x}$ is defined on at most one element of each orbit of $U(S^{Ev})$. Namely, if $i \in Dom(x), x \in S$, then the unique element of the orbit $g\mathcal{O}_{i}$ on which $h_{x}$ is defined is $(g,(i,i))$ and is not defined on $g\mathcal{O}_{i}$ if $i$ is not in $\operatorname{Dom}(x)$.} 

\item{$\operatorname{Dom}(h_{x}) \cap Im(h_{x})=\emptyset$ for all $x \in S$.}\label{p3}

\item{We make the previous point more precise. For $k,l = 0, \ldots n$, let $B_{k}=\{(i,i+k)|i=1, \ldots n\}$ and $R_{l}=\{(i,l)|i = 1, \ldots n\}$. Then:

\medskip

For all $x\in S$ we have $\operatorname{Dom}(h_{x}) \subseteq G \times B_{0}, Im(h_{x}) \subseteq G \times R_{0}$ and point \ref{p3}. follows.}

\item{$B_{k}t=B_{k+1}, R_{l}t=R_{l+1}, k,l = 0 \ldots n$. }

\item{The collections $\{G \times B_{k}|k = 0 \ldots n\}$ and $\{G \times R_{l} | l=0\ldots n\}$ are blocks (in the sense of the theory of permutation groups- equivalently partition classes of right congruences) of the permutation group $(G \times B^{Ev}, Z_{n+1})$.}

\end{enumerate}

\bds

Before continuing we discuss the intuition behind the construction of $S^{Ev}$. We have replaced the partial functions in $S$ acting on $G \times B$ by the collection of partial functions $\{h_{x} \mid x \in S\}$ acting on $G \times B^{Ev}$. Each $h_{x}$ has domain (as a monomial function) contained in the ``diagonal" set $B_{0}=\{(i,i)\mid i=1\ldots n\}$ and image contained in $R_{0}=\{(i,0)\mid i=1\ldots n\}$. Clearly all of these functions compose to the empty function. We have remembered the formula for each $x \in S$ and stored the values in $G \times R_{0}$ but forgotten the composition in $S$. $t$ acts as a $\pmod{n+1}$ clock, mapping $R_{i}$ to $R_{i+1}$. Note then that $R_{i}, i>0$ intersects the diagonal in precisely one point $(i,i)$. It follows that any element of $S^{Ev}$ that as a composition of generators uses at least 2 of the $h_{x}$ has rank 1 as a monomial function. This will imply that $S^{Ev}$ is a smallish monoid. We then will see that the operator $\omega+*$ in the evaluation semigroup $Eval(S^{Ev})$ applied to $t$ allows us to ``recover" the formula and composition rules for elements of $S$. We give the details below.

We see that the main idea of our construction is to replace a collection $X$ of partial functions from a set $Q$ to itself by a collection $h_{X}$ of partial functions between disjoint sets $D_{Q}$ to $R_{Q}$ in bijection with $Q$. If $f\colon Q \rightarrow Q \in X$ then we define $h_{f}:D_{Q} \rightarrow R_{Q} \in h_{X}$ by $d_{q}h_{f}=r_{qf}$ for $q \in \dom(f)$. 

The idea of separating the domain and range of functions occurs in a number of fields of mathematics.
In graph theory this leads to the bipartite double of a graph. In category theory it is a special case of the cograph (or collage) of a profunctor. The reader can search the literature for these and other cases of this construction.

As mentioned above, for us, this allows us to store the values of the elements of $S$ in $G \times R_{0}$ and then use $t$ to "mutate" these values in order to build a smallish monoid. 
\eds

\bl

$ S^{Ev}$ acts transitively on $G \times B^{Ev}$.

\el

\begin{proof}

Let $(g,(i,j)), (h,(k,l)) \in G \times B^{Ev}$. Choose $a' \in A$ such that $C(i,a') \neq 0$ and let $x =(a',C(i,a')^{-1}g^{-1}h,k) \in I(S)$. Then $(g,(i,j))t^{j-i}h_{x}=(g,(i,i)h_{x})=(h,(k,0))$. Therefore, \\ $(g,(i,j))t^{j-i}h_{x}t^{l} = (h,(k,l))$. 

\end{proof}

It is well known \cite{Arbib, qtheory}  that since $S^{Ev}$ acts faithfully and transitively on the right of a set, $S^{Ev}$ has a unique 0-minimal ideal $I(S^{Ev})$ that is a 0-simple subsemigroup of $S^{Ev}$.
$I(S^{Ev})$ consists of all elements in $S^{Ev}$ whose image on $G \times B^{Ev}$ is either the empty set or $G \times (i,j)$ for some $(i,j) \in B^{Ev}$. We will determine the structure of $I(S^{Ev})$ below. We first prove that $S^{Ev}$ is a smallish monoid.

\bl \label{smallish}

Let $I$ be the ideal of $S^{Ev}$ generated by $\{h_{x} \mid x \in S\}$.

\begin{enumerate}

\item{$I$ is the ideal of non-units $S^{Ev} - U(S^{Ev})$ of $S^{Ev}$.}
    
\item{ $I^{2} = I(S^{Ev})$.}

\end{enumerate}

\el

\begin{proof}

Item 1. is clear since $t$ generates the group of units of $S^{Ev}$ and all $h_{x}$ are non-invertible having rank at most $|G|n$ whereas $|G \times B^{Ev}| = |G|n(n+1)$.

Since $I(S^{Ev})$ is a regular semigroup, it follows that $I(S^{Ev})$ is contained in every power of $I$ and in particular, $I(S^{Ev}) \subseteq I^2$. For the opposite inclusion it's enough to show that elements of the form $t^{k}h_{x}t^{l}h_{y}, k,l\geq 0, x,y \in S$ belong to $I(S^{Ev})$. 

Let $(i,j) \in B^{Ev}$. Then $(i,j)t^{k}h_{x} = (i,j+k)h_{x} = g_{i,x}(ix,0)$ if $i=j+k$ and $i \in Dom(x)$ and undefined otherwise. Therefore, $D = Dom(t^{k}h_{x}t^{l})= G \times \{(i,i-k) \mid i \in Dom(x)\}$. It follows that if $(i,i-k) \in D$, then 
$(i,i-k)t^{k}h_{x}t^{l}=g_{i,x}(ix,l)$. But $(ix,l) \in Dom(h_{y})$ if and 
only if $ix = l$ and $l \in Dom(y)$. Since $l$ is fixed, we have that 
$t^{k}h_{x}t^{l}h_{y}$ is the empty function if $l \not\in Dom(y)$ and has 
image $G \times \{(ly,0)\}$ otherwise. Therefore, either
$\Im(t^{k}h_{x}t^{l}h_{y})$ is empty or $\Im(t^{k}h_{x}t^{l}h_{y})= G \times \{(ly,0)\}$ and thus belongs to the 0-minimal ideal of $S^{Ev}$.

\end{proof}

\bc

$S^{Ev}$ is a smallish monoid.

\ec

\begin{proof}

This follows immediately from Lemma \ref{smallish} and the definition of smallish monoid.


\end{proof}

We wish to describe $I(S^{Ev})$ as a Rees matrix semigroup. We first claim that the maximal subgroup of $I(S^{Ev})$ is $G$. To prove this, let $x = (a,g,1)$ be an element of the distinguished $\mathcal{H}$-class in $I(S)$. Recall that $C(1,a)=1$ and that $C(i,a) \in \{0,1\}$ for $i=1, \ldots , n$. Then $h_{x}t$ sends $(i,i)$ to $(g,1)$ if $C(i,a)=1$ and undefined if $C(i,a)=0$. Therefore, $\Im(h_{x}t)=G \times\{(1,1)\}$ and thus belongs to $I(S^{Ev})$. It follows from $C(1,a)=1$ that the map sending $g$ to $h_{(a,g,1)}t$ is an isomorphism to an $\mathcal{H}$-class of the 0-minimal ideal of $RM_{B^{Ev}}(G)$.  $RM_{B^{Ev}}(G)$ has a unique minimal 0-ideal  $I(RM_{B^{Ev}}(G))$ with maximal subgroup $G$ and whose non-zero elements have image $G \times (i,j)$ for some $(i,j) \in B^{Ev}$. Therefore, since  $I(S^{Ev}) \subseteq I(RM_{B^{Ev}}(G))$, it follows that the maximal subgroup of $I(S^{Ev})$ is $G$.

Since $S^{Ev}$ acts transitively on $G \times B^{Ev}$ and the maximal subgroup of $I(S^{Ev})$ is $G$ it follows that $I(S^{Ev})$ is isomorphic to a Rees matrix semigroup of the form $\mathcal{M}^{0}(G,A^{Ev},B^{Ev},C^{Ev})$. We first determine $A^{Ev}$. 

$A^{Ev}$ indexes the $\mathcal{R}$-classes of $I(S^{Ev})$. Let $s \in I(S^{Ev})$. Then $\Im(s) = G \times \{(i,j)\}$ for some $(i,j) \in B^{Ev}$. The $\mathcal{R}$-class of $s \in I(S^{Ev})$ is determined by the kernel of $s$ denoted by $\operatorname{ker}(s)$, where we are treating $s$ as a partial function on $G \times B^{Ev}$. Recall that $\operatorname{ker}(s)$ is the equivalence relation on $\operatorname{Dom}(s)$ such that $((h,(k,l)),(k,(u,v))) \in ker(s)$ if and only if $(h,(k,l))s=(k,(u,v))s$. Equivalently, the partition classes of $\operatorname{ker}(s)$ are given by the fibers $(g,(i,j))s^{-1}$ where $(g,(i,j)) \in Im(s)$. Note that 
$(g,(i,j))s^{-1}=g((1,(i,j))s^{-1})$. We consider $(1,(i,j))s^{-1}$ as a partial function $f_{i,j,s}:B^{Ev}\rightarrow G$, where $(k,l)f_{i,j,s}=g_{k,l}$ if $(g_{k,l},(k,l)) \in (1,(i,j))s^{-1}$ and undefined otherwise. Therefore the cross-section in the sense of the Rhodes lattice \cite{AmigoDowling} defined by $f_{i,j,s}$ determines the kernel of $s$. 

Thus we need to determine which cross-sections occur as $\mathcal{R}$-classes of $S^{Ev}$. The next Proposition does precisely this. For conciseness and as usual we write $(i,j)$ for $(1,(i,j)) \in G \times B^{Ev}$.

\bp \label{RofSEv}

Let $(G \times B,S)$ be a $\GM$ transformation semigroup with $I(S) = \mathcal{M}^{0}(G,A,B,C)$. Then the following holds.

\begin{enumerate}\label{ISEv}

\item{Let $x = (a',g,j) \in I(S)$ and $k,l \in \{0 , \ldots n\}$. Let $s = t^{k}h_{x}t^{l} \in S^{Ev}$. Then $\operatorname{Dom}(s)=G \times \{(i,i-k) \mid C(i,a') \neq 0\}$. If $(h,(i,i-k)) \in Dom(s)$, then $(h,(i,i-k))s = (hC(i,a')g,(j,l))$. Thus the image of $s$ acting on 
$G \times B^{Ev}$ is $G \times \{(j,l)\}$. Therefore $s \in I(S^{Ev})$.  The domain of the cross-section $f_{j,l,s}$ is given by $\{(i,i-k) \mid C(i,a') \neq 0\}$ and for $(i,i-k) \in Dom(s)$ we have $((i,i-k))f_{j,l,s}=C(i,a')g$.} 

\item{Let $x,y \in S$ and let $k,l \in \{0 , \ldots n\}$. Let $s =t^{k}h_{x}t^{l}h_{y}$. Then $s$ is not the empty function if and only if
$ix = l$ and $l \in Dom(y)$. In this case, $\operatorname{Dom}(s)=G \times \{(i,i-k) \mid i \in lx^{-1}\}$. If $(h,(i,i-k)) \in Dom(s)$, then 
$(h,(i,i-k))s = hg_{i,x}g_{l,y}(ly,0)$. Thus the range of $s$ is $G \times \{(ly,0)\}$ and $s \in I(S^{Ev})$. The cross-section defined by $s$ sends $(i,i-k)$ to $g_{i,x}g_{l,y}$. 
  }

\item{Every element of $I(S^{Ev})$ is $\mathcal{R}$-equivalent to an element of one of the forms in Item 1. or Item 2.}

\end{enumerate}

\ep

\begin{proof}

Item 1. and Item 2. can be proved by direct calculation. To prove Item 3., let $s \in I(S^{Ev})$. Then as a product of generators, $s$ must contain at least 1 generator of the form $h_{x}$, where $x \in S$. Otherwise, $s$ is a power of $t$ and is in the group of units of $S^{Ev}$. If $s =t^{k}h_{x}t^{l} \in I(S^{Ev})$ and is not the empty function, we claim that $x \in I(S)$. Assume that $x$ is not in $I(S)$ then since $t$ is a unit, $\operatorname{rank}(s) = rank(h_{x})=rank(x) > |G|$. Since the non-zero  of $I(S^{Ev})$ have rank $|G|$, it follows that $x \in I(S)$ and $s$ is as in Item 1.

If $s$ is a product of generators that has at least two occurrences of generators other then $t$, then $s=t^{k}h_{x}t^{l}h_{y}v$ for some 
$x,y \in S, v \in S^{Ev}$. By Lemma \ref{smallish}, $t^{k}h_{x}t^{l}h_{y} \in I(S^{Ev})$. Therefore if $s \neq 0$, then
$s \mathcal{R} t^{k}h_{x}t^{l}h_{y}$ since $I(S^{Ev})$ is a 0-simple semigroup. This completes the proof.

\end{proof}

\bc

Let $A^{Ev}$ be the set of cross-sections of all elements of the form described in Item 1. and Item 2. of Proposition \ref{ISEv}. Then
$A^{Ev}$ is an index of the set of $\mathcal{R}$-classes of $I(S^{Ev})$.

\ec

The next Lemma gives the connection between an element in $I(S^{Ev})-\{0\}$ and the corresponding cross-section describing its $\mathcal{R}$-class. We state it for a general $\GM$ semigroup and will apply it in describing the structure matrix of $S^{Ev}$.

\bl\label{Rclass}

Let $(G \times B,S)$ be a $\GM$ transformation semigroup with $I(S) = \mathcal{M}^{0}(G,A,B,C)$ and  $B =\{1, \ldots n\}$. 

\begin{enumerate}
    
\item{Let $s \in I(S)$ have image $\Im(s) =G \times \{j\}$ and defined by the row-monomial matrix that sends $i$ to $g_{i}j$ for $i \in Dom(s)$. Then $(1,j)s^{-1} = \{(g_{i}^{-1},i)\mid i \in Dom(s)\}$.} 

\item{Considering $(1,j)s^{-1}$ as the function $f_{s}:Dom(s) \rightarrow G$ with $if_{s}=g_{i}^{-1}$, the $\mathcal{R}$-class of $s$ is the set of all row-monomial matrices $t \in I(S)$ that send $i \in Dom(s)$ to $(g(if_{s}))^{-1}b = (g_{i}g^{-1})b$, where $g \in G, b \in B$. This is precisely the set of elements $t$ of $I(S)$ such that the cross-section defined by $bt^{-1}$ is $[f_{s}]$.}

\end{enumerate}
\el
 \begin{proof}

Both items are proved by a straightforward calculation.
     
 \end{proof}
We now describe the structure matrix of $I(S^{Ev})$. Let $k \in \{0,\ldots ,n\}$. Recall that $B_{k} \subseteq B^{Ev}$ is defined  by $B_{k}=\{(i,i+k)\mid i = 1, \ldots n\}$ where $i+k$ is taken modulo $n+1$. Then $B^{Ev}$ is the disjoint union of all the $B_{k}$. Furthermore it follows from Proposition \ref{RofSEv} that the domain of every element $s \in I(S^{Ev})$ is contained in $G \times B_{k}$ for a unique $k$. That is, the $B_{k}$ partition $B^{Ev}$ and also the domains of elements in $I(S^{Ev})$.

Let $A_{k}$ be the set of cross-sections of elements of $I(S^{Ev})$ whose domain is contained in $B_{k}$. We choose a representative $a$ for each cross-section $[a]\in A^{Ev}$. For each $b \in B^{Ev}$ there is a unique $s \in I(S)$ whose image is $b$ and such that $bs^{-1} = a$. It follows that $I(S^{Ev})$ is the disjoint union of $A_{k} \times G \times B_{l}$ where $1 \leq k,l \leq n$.

Suppose that $k \neq l$ and let $s \in A_{k} \times G \times B_{l}$. Then the unique element in the range, $(i,i+l)$ is not in the domain of $s$ since the domain of $s$ is contained in $B_{k}$ and $k \neq l$. Therefore $s^{2}=0$. It follows that no $\mathcal{H}$-class in $A_{k} \times B_{l}$ is a group. This means that all the entries of the structure matrix in places corresponding to $B_{l} \times A_{k}$ are 0. Therefore, the structure matrix is the direct sum of the structure matrices $C_{k}^{Ev}$ corresponding to the submatrix in positions $B_{k} \times A_{k}$.

For each $k$,  $A_{k} \times G \times B_{k} \cup \{0\}$ is a subsemigroup of 
$I(S^{Ev})$.  We claim that 
$A_{k} \times G \times B_{k} \cup \{0\} = t^{-k}(A_{0} \times G \times B_{0})t^{k} \cup \{0\}$. Indeed if $s \in (A_{0} \times G \times B_{0})$ and $(i,i)s=(j,j)$, then $(i,i+k)t^{-k}st^{k}=(j,j+k)$. Therefore all the $(A_{k} \times G \times B_{k}) \cup \{0\}$ are isomorphic to each other. It suffices to then construct the structure matrix corresponding to positions in $B_{0} \times A_{0}$.

We claim that that $A_{0} \times G \times B_{0} \cup \{0\}$ is a regular subsemigroup of $S^{Ev}$. Let $[a] \in A_{0}$ and let $(i,i)$ be an element in the domain of $a$. Then the function sending all elements in $a$ to $(i,i)$ defines an idempotent. So each $\mathcal{R}$-class of $A_{0} \times G \times B_{0} \cup \{0\}$ contains an idempotent. Now let $e=(a',C(i,a')^{-1},i)$ be an idempotent of $S$. Then $h_{e}t^{i}$ is an idempotent in $A_{0} \times G \times B_{0} \cup \{0\}$ with image $(i,i)$. Therefore each $\mathcal{L}$-class of $A_{k} \times G \times B_{k} \cup \{0\}$ contains an idempotent. This proves the claim.

We choose Rees representatives for the semigroup $A_{0} \times G \times B_{0}\cup \{0\}$. Recall that $B_{0} =\{(i,i)\mid i=1, \ldots, n\}$ and we order $B_{0}$ with $(1,1) < (2,2) \ldots < (n,n)$. Let $b_{i}=(a,1,i)$ be the Rees representative in $S$ for the $\mathcal{L}$-class with index $i$ in $I(S)$. We take $h_{b_{i}}t^{i}$ a representative for the $\mathcal{L}$-class $(i,i)$.  As a row-monomial function, $h_{b_{i}}t^{i}$ has domain $\{(j,j)\mid C(j,a) \neq 0\}$ and such a $(j,j)$ is sent to $C(j,a)(i,i)$.

Let a be a cross-section $[a] \in A_{0}$ as described in Proposition \ref{RofSEv}. We choose a representative $a$ of the cross-section such that the minimal element in its domain (in the order on $B_{0}$ above) has value 1. We now take as Rees representative of the $\mathcal{R}$-class indexed by $[a]$, the element that sends $(j,j)$ to $((j,j)a)^{-1}(1,1)$ for all $(j,j) \in Dom(a)$. 

By the proof of the Rees Theorem the structure matrix $C_{0}$ is defined by $C_{0}((i,i),a)= g$  if $(i,i) \in Dom(a)$ and $(1,1)h_{b_{i}}t^{i}a=g(1,1)$ and 0 otherwise. We summarize this discussion in the next Proposition.

\bp\label{structmat}
The structure matrix $C^{Ev}$ of $I(S^{Ev})$ is the direct sum of $n+1$ copies of the matrix $C_{0}$ defined above.
\ep

As an example, let $s= (a',1,i) \in I(S)$. Then $is^{-1}=\{(C(j,a')^{-1},j) \mid C(j,a') \neq 0\}$. Therefore, the cross-section associated to $s$ is defined by $jf_{i,s}=C(j,a')^{-1}$. By Proposition \ref{RofSEv} the row 
of $C_{0}$ corresponding to the $\mathcal{R}$-class of $h_{s}t^{i}$ is the same as row $a'$ of $C$. Thus $C$ is the submatrix of $C_{0}$ consisting of the columns corresponding to elements in $I(S)$ when we identify $i$ with $(i,i)$. In other words, $C_{0}$ is obtained from $C$ by adding columns corresponding to cross-sections defined by Proposition \ref{RofSEv}. 

Structure matrices are unique only up to permutation of rows and columns and multiplying rows on the left and columns on the right by elements of $G$ (this is in fact cohomological equivalence of the corresponding Graham-Houghton graph \cite[Section 4.13]{qtheory}). In order to see  the direct sum structure given by Proposition \ref{structmat} we have to order $B^{Ev}$ and $A^{Ev}$ in an appropriate way. There is a natural order for $B^{Ev}$ that we now define. Unfortunately, there is not a natural way to order $A^{Ev}$.

We order $B^{Ev}$ as follows. For $k= 0 , \ldots , n$ we have defined 
$B_{k}=\{(i,i+k) \pmod{n+1}\} \mid i =  1, \ldots n\}$. We order $B_{k}$ by $(i,i+k) \leq (j,j+k)$ if $i \leq j$. We then order $B^{Ev}$ so that $B_{0} < B_{1} \ldots < B_{n}$. More concretely the order on $B^{Ev}$ is:

$$(1,1), (2,2), \ldots (n,n), (1,2), (2,3), \ldots , (n,0) , (1,3) \ldots , (1,0), \ldots , (n,n-1)$$ 

We note that for $l=0 \ldots n$, $B_{0}t^{l}=B_{l}$. We have already chosen the Rees representatives $h_{b_{i}}t^{i}, i = 1, \ldots, n$ for $B_{0}$. We can therefore take $\{h_{b_{i}}t^{i+l}\}, i = 1, \ldots, n$ as Rees representatives for $B_l$. Therefore $\{h_{b_{i}}t^{i+l} \mid i = 1, \ldots , n, l=0 \ldots , n\}$ is a system of Rees representatives for $B^{Ev}$. We fix an order for the Rees representatives we picked for $A_{0}$. We note that for each $k = 0, \ldots n$ , $t^{-k}A_{0}=A_{k}$ and we take $t^{-k}A_{0}$ as Rees representatives for $A_{k}$. 
The proof we gave showing that  
$A_{k} \times B_{k} = t^{-k}(A_{0} \times B_{0})t^{k}$ gives the desired direct sum decomposition from Proposition \ref{structmat}. See Figure \ref{MatrPict} for a picture of the structure matrix.

\setcounter{figure}{0}
\begin{figure}[ht]
\centering
\begin{tikzpicture}[baseline]

  \matrix (M) [matrix of math nodes,
              nodes in empty cells,
              column sep=1.0cm,
              row sep=0.8cm,
              anchor=west] at (0.2,0.0)
  {
    C_0 & 0 & \cdots & 0 & \cdots & 0 \\
    0 & C_0 & \cdots & 0 & \cdots & 0 \\
    \vdots & \vdots & \ddots & \vdots & & \vdots \\
    0 & 0 & \cdots & C_0 & \cdots & 0 \\
    \vdots & \vdots & & \vdots & \ddots & \vdots \\
    0 & 0 & \cdots & 0 & \cdots & C_0 \\
  };
  \node[anchor=south] at ($(M-1-1.north)+(0,0.35)$) {$A_0$};
  \node[anchor=south] at ($(M-1-2.north)+(0,0.35)$) {$A_1$};
  \node[anchor=south] at ($(M-1-3.north)+(0,0.35)$) {$\cdots$};
  \node[anchor=south] at ($(M-1-4.north)+(0,0.35)$) {$A_i$};
  \node[anchor=south] at ($(M-1-5.north)+(0,0.35)$) {$\cdots$};
  \node[anchor=south] at ($(M-1-6.north)+(0,0.35)$) {$A_n$};

  \node[anchor=east] at ($(M-1-1.west)+(-0.35,0)$) {$B_0$};
  \node[anchor=east] at ($(M-2-1.west)+(-0.35,0)$) {$B_1$};
  \node[anchor=east] at ($(M-3-1.west)+(-0.35,0)$) {$\vdots$};
  \node[anchor=east] at ($(M-4-1.west)+(-0.35,0)$) {$B_i$};
  \node[anchor=east] at ($(M-5-1.west)+(-0.35,0)$) {$\vdots$};
  \node[anchor=east] at ($(M-6-1.west)+(-0.35,0)$) {$B_n$};

  \coordinate (NW) at ($(M.north west)+(-0.25,0.25)$);
  \coordinate (SW) at ($(M.south west)+(-0.25,-0.25)$);
  \coordinate (NE) at ($(M.north east)+(0.25,0.25)$);
  \coordinate (SE) at ($(M.south east)+(0.25,-0.25)$);

  \draw[line width=0.6pt] (NW) -- (SW);
  \draw[line width=0.6pt] (NW) -- ($(M.north west)+(0,0.25)$);
  \draw[line width=0.6pt] (SW) -- ($(M.south west)+(0,-0.25)$);

  \draw[line width=0.6pt] (NE) -- (SE);
  \draw[line width=0.6pt] (NE) -- ($(M.north east)+(0,0.25)$);
  \draw[line width=0.6pt] (SE) -- ($(M.south east)+(0,-0.25)$);

  \node[anchor=north] at ($($(M-6-1.south west)!0.5!(M-6-6.south east)$)+(0,-0.8)$)
    {$B_i=B_0t^{i},\;\; A_i=t^{-i}A_0,\; i=0,\ldots,n$};
\end{tikzpicture}
\caption{The structure matrix of $I(S^{Ev})$.}\label{MatrPict}
\end{figure}

We now use the description of the structure matrix to prove that the small submonoid $Sm(S^{Ev})$ always has complexity 1. Recall that $Sm(S^{Ev})$ is the union of the group of units $U(S^{Ev})$ and the unique 0-minimal ideal 
$I(S^{Ev})$. We first recall Tilson's 2-$\mathcal{J}$-class Theorem that determines the complexity of a small monoid.
In fact, the theorem gives a decidability criterion for all semigroups with 2 non-zero $\mathcal{J}$-classes. The general case can be reduced to that of small monoids. See \cite{2J}, 
\cite[Section 4.15]{qtheory} for details. From the depth decomposition theorem \cite{TilsonXI}, a small monoid has complexity at most 2. Thus we need a criterion to decide whether the complexity is 1 or not.

Let $S = H \cup M^{0}(G,A,B,C)$ be a small monoid. We identify $B$ with the set of $\mathcal{L}$-classes of $M^{0}(G,A,B,C)$. Then $H$ acts on the right of $B$ by $L\cdot h = Lh$, where $L$ is an $\mathcal{L}$-class of $M^{0}(G,A,B,C)$. This is the restriction of the action defining $\RLM(S)$ to $H$. Then for each such $L$, $H \cup LH \cup \{0\}$ is a submonoid of $M$ called a right-orbit monoid. The collection of all $LH$ are the orbits of the action of $H$ on $B$. Let $k=k(S)$ be the number of such orbits. Let $B_{1}, \ldots B_{k}$ be the orbits of $H$ on $B$. Then $G \times B_{1}, \ldots , G \times B_{k}$ is a partition of $G \times B$. Furthermore the right-orbit monoids are exactly the monoids $H \cup A \times G \times B_{i}, i=1, \ldots k$. Notice that $A \times G \times B_{i} = M^{0}(G,A,B_{i},C_{i})$, where $C_{i}$ is the restriction of $C$ to $B_{i} \times A$. This is a not necessarily regular Rees matrix semigroup.

\bt\label{2J}

Let $S = H \cup M^{0}(G,A,B,C)$ be a small $\GM$ monoid. Then $Sc = 1$ if and only if each right-orbit monoid $H \cup LH \cup \{0\}$ has aperiodic idempotent generated subsemigroup. 

\et

We compute orbit monoids in $Sm(S^{Ev})$. An $\mathcal{L}$-class is indexed by a unique $(i,j) \in B^{Ev}$. The orbit $(i,j)Z_{n+1}=\mathcal{O}_{i}=\{(i,k)\mid k=0,\ldots n\}$. Using the notation above (see Figure \ref{MatrPict}) $(i,k) \in B_{k-i}$. Thus the orbit $\mathcal{O}_{i}$ intersects exactly 1 row of each $B_{l}, l=0,\ldots n$ in the direct sum decomposition in Figure \ref{MatrPict}. Thus the orbit is a disjoint union of $\mathcal{L}$-classes, 
$\{L_{i}\mid i=0,\ldots n\}$ such that for $i\neq j$, it follows from the direct sum decomposition that $L_{i}L_{j}=\{0\}$. It follows from Rees Theorem that the idempotents in a regular $\mathcal{L}$-class of any semigroup form a left-zero semigroup. Therefore the idempotent generated submonoid of the orbit monoid consists of the identity and a 0-disjoint union of left-zero semigroups and is thus aperiodic. This proves the following theorem.

\bt \label{small1}

The small submonoid $Sm(S^{Ev})$ has complexity 1.

\et

\begin{proof}

Since the minimal ideal $I(S^{Ev})$ contains the non-trivial group $G$, $Sm(S^{Ev})c >0$. The discussion preceding the statement of the theorem shows that the idempotent generated submonoid of every orbit monoid is aperiodic and the result follows from Theorem \ref{2J}.
    
\end{proof}
The following Proposition improves Proposition 5.16 of \cite{Trans}. The reader should review the definition of the operator $\omega+*$ and the evaluation transformation semigroup of a $\GM$-transformation semigroup. See Section 4 of \cite{Trans} or Appendix \ref{Eval} of this paper.

\bp\label{EmbedEval}

Let $X=(G \times B,S)$ be a $\GM$-transformation semigroup with $B=\{1,\ldots n\}$. Then $X$ embeds into the evaluation transformation semigroup $\mathcal{E}(S^{Ev}) = (States(S^{Ev}), Eval(S^{Ev})).$ 

\ep

\begin{proof}

For each $g\in G$ we have $(g,(i,0))t^{\omega+*}$ is the orbit $g\mathcal{O}_{i}=\{(g,(i,j))\mid j=0,\ldots n\}$ . Thus every orbit $g\mathcal{O}_{i} \in States(Eval(S^{Ev}))$ where the cross-section takes the constant value 1. Let $x \in S$. Then for each $i \in Dom(x)$ we have $(g,(i,0))t^{\omega+*}h_{x}=(gg_{i,x},(ix,0))$. Let $S'$ be the subsemigroup of $Eval(S^{Ev})$ generated by $t^{\omega+*}h_{x}, x \in S$ and let $B'= \{(i,0)|i = 1 , \ldots n\}$. It follows that $(G \times B',S')$ is a subtransformation semigroup of $(States(S^{Ev}), Eval(S^{Ev}))$ that is isomorphic to $(G \times B,S)$. 

\end{proof}

We now arrive at the main result of this section. This will prove all of Theorem \ref{evalembed}.

\bt

A $(G \times B,S)$  $\GM$ transformation semigroup has an aperiodic flow if and only if $(G \times B^{Ev},S^{Ev})$ has an aperiodic flow.

\et

\begin{proof}

Assume that $S^{Ev}$ has an aperiodic flow. Then by  \cite[Section 4]{Trans}, no state of $\mathcal{E}(S^{Ev}) = (States(S^{Ev}), Eval(S^{Ev}))$ is the contradiction of the Rhodes lattice. 
Let $x \in S$. Then for each $i \in Dom(x)$ we have $(i,0)t^{\omega+*}h_{x}=(ix)(ix,0)$. Let $S'$ be the subsemigroup of $Eval(S^{Ev})$ generated by $t^{\omega+*}h_{x}, x \in S$ and let $B'= \{(i,0)|i = 1 , \ldots n\}$. Proposition \ref{EmbedEval} states that $(G \times B',S')$ is a subtransformation semigroup of $(\operatorname{States}(S^{Ev}), \operatorname{Eval}(S^{Ev}Ev))$ that is isomorphic to $(G \times B,S)$. Therefore, every state of $\mathcal{E}(S)= (States(S), Eval(S))$ also belongs to $\operatorname{States}(S^{Ev})$ and is therefore not the contradiction. It follows from Theorem 3.4.1 of \cite{complexity1} that $(G \times B,S)$ has an aperiodic flow.

Now assume that $(G \times B,S)$ has an aperiodic flow. Let $(Q,T)$ be an aperiodic transformation semigroup and $F:Q \rightarrow Rh_{B}(G)$
be a flow. Without loss of generality we will assume that $T$ is a monoid. We will now extend $F$ to a flow 
$F^{Ev}:Q \rightarrow Rh_{B^{Ev}}(G)$.

For $i \in B$, let $i^{Ev}=\{(i,j) \mid j = 0, \ldots n\}$ be the orbit of $\la t \ra$ containing $(i,0)$. Note that 
$i^{Ev}=(i,0)t^{\omega+*}$. If $X \subseteq B$, let $X^{Ev}=\bigcup_{i \in X}i^{Ev}$ be the union of the orbits containing the 
elements of $X$. If $\Pi$ is a partition of $X \subseteq B$ with partition classes $\{\pi_{1}, \ldots \pi_k\}$, let 
$\Pi^{Ev}= \{\pi_{1}^{Ev}, \ldots \pi_k^{Ev}\}$. If $\Theta=(X,\Pi,[f]) \in Rh_{B}(G)$, let 
$\Theta^{Ev}=(X^{Ev},\Pi^{Ev},[f^{Ev}]) \in Rh_{B^{Ev}}(G)$ where we define $[f^{Ev}]$ as follows. Let $\pi_{j}^{Ev}$ be a class of $\Pi^{Ev}$.
Let $(k,l) \in \pi_{j}^{Ev}$. Then we let $f^{Ev}((k,l))=f_{j}(k)$, where $f_{j}$ is the component of $f$ defined on $\pi_{j}$.

We define a function $Q \rightarrow Rh_{B^{Ev}}(G)$ by $qF^{Ev}=(qF)^{Ev}$. We verify that $F^{Ev}$ is a flow. Let $qF^{Ev}= (X^{Ev},\Pi^{Ev},[f^{Ev}])$ where $qF=(X,\Pi,[f])$. We have $X^{Ev}t=X^{Ev}$, since $X$ is a union of orbits of $t$ by definition. Similarly, $t$ fixes every partition class of $\Pi^{Ev}$ and $t$ therefore preserves cross-sections as well. So $t$ is covered by the identity element of $T$.

Now let $x \in S$. Let $x$ be covered by $\widehat{x} \in T$. Let $(q\widehat{x})F = (Y,\Theta,[k])$. Then multiplication by $x$ induces a 
one-to-one partial map $(X,\Pi,[f])$ to $(Y,\Theta,[k])$ that maps the cross-section $[f]$ into the cross-section $[k]$. Now consider $h_{x}$ a generator of $S^{Ev}$. We cover $h_{x}$ by $\widehat{x}$. Since $h_{x}$ is defined on at most one element of each orbit of $t$, it is immediate from the definition that right multiplication by $h_{x}$ induces an injective function from $(X^{Ev},\Pi^{Ev},[f^{Ev}])$ to $(Y^{Ev},\Theta^{Ev},[k]^{Ev})$ that maps the cross-section $[f]^{Ev}$ into the cross-section $[k]^{Ev}$. Therefore $F^{Ev}$ is a flow and the theorem is proved.

\end{proof}

\bc \label{reduction}

The problem of deciding if a $\GM$ transformation semigroup $(G \times B,S)$ has an aperiodic flow can be reduced to the case that $S$ is a smallish monoid.

\ec
 
\bc

Let $S$ be a $\GM$ semigroup such that $\RLM(S)c \leq 1$. Then $Sc = 1$ if and only if $S^{Ev}c=1$. Thus the problem of deciding if a $\GM$ semigroup has complexity 1 is reducible by induction to checking that $\RLM(S)c$ has complexity at most 1 and whether the problem of deciding if a smallish $\GM$ monoid has complexity 1.

\ec

\begin{proof}

Since as for any $\GM$ smallish monoid we have $\RLM(S^{Ev})c = 1$. Therefore, by the main theorem of \cite{complexity1}, $S^{Ev}c=1$ if and only if $S^{Ev}$ has an aperiodic flow. The result now follows from Corollary \ref{reduction}.
 
\end{proof}
\begin{center}
{\LARGE Appendices}\label{Append}
\end{center}

\appendix

\section{APPENDIX: Flows and the Flow Decomposition Theorem}\label{Flows}

We review the definition of flows from  an automaton with alphabet $X$  to the set-partition $SP(G\times B)$ and to the Rhodes lattice $Rh_{B}(G)$. For more details, see \cite[Sections 2-3]{Trans}. We first need to recall details about the Type II subsemigroup of a semigroup.

\subsection{The Type II Subsemigroup and the Tilson Congruence}\label{redux.sec}

In this subsection we review the type II subsemigroup $S_{II}$ of a semigroup. It plays an important role in finite semigroup theory and a central role in flow theory. It was first defined in \cite{lowerbounds2}. In that paper it was proved that it is decidable if a regular element of a semigroup $S$ belongs to $S_{II}$. In particular, it followed that if $S$ is a regular semigroup, then membership in $S_{II}$ is decidable. Later 
Ash \cite{Ash} proved that membership in $S_{II}$ is decidable for all finite semigroups $S$.

The type II subsemigroup $S_{II}$ of $S$ is the smallest subsemigroup of $S$ containing all idempotents and closed under weak conjugation: if $xyx = x$ for $x,y \in S$, then $xS_{II}y \cup yS_{II}x \subseteq S_{II}$. Membership in $S_{II}$ is clearly decidable from this definition. Its importance stems from the following Theorem, where we give its original, not a priori decidable definition.

\bt

Let $S$ be a finite semigroup. Then $S_{II}$ is the intersection of all subsemigroups of the form $1\phi^{-1}$ where $\phi:S \rightarrow G$
is a relational morphism from $S$ to a group $G$.

\et

The aforementioned decidability results in \cite{lowerbounds2, Ash} prove that the two definitions we give define the same subsemigroup of a semigroup $S$. Here are some important properties of the type II subsemigroup. See \cite{qtheory} for proofs.

\bt\label{typeIIprop}

\begin{enumerate}

\item{Let $f:S \rightarrow T$ be a morphism between semigroups $S$ and $T$. If $s \in S_{II}$, then $sf \in T_{II}$. Therefore the restriction of $f$ to $S_{II}$ defines a morphism $f_{II}:S_{II} \rightarrow T_{II}$. This defines a functor on the category of finite semigroups.}

\item{If $S$ divides $T$, then $S_{II}$ divides $T_{II}$.}

\item{Assume that a semigroup $S$ divides a semidirect product $T*G$ where $T$ is a semigroup and $G$ is a group. Then $S_{II}$ divides $T$.}

\end{enumerate}

\et





A {\em congruence} on a transformation semigroup $(Q,S)$ is an equivalence relation $\approx$ on $Q$ such that if $q\approx q'$ and for $s \in S$, and both $qs$ and $q's$ are defined, then $qs\approx q's$. Every $s \in S$ defines a partial function on $\faktor{Q}{\approx}$ by 
$[q]_{\approx}s=[q's]_{\approx}$ if  $q's$ is defined for some $q' \in [q]_{\approx}$. The quotient $\faktor{(Q,S)}{\approx}$ has states $\faktor{Q}{\approx}$ and semigroup $T$ the semigroup generated by the action of all $s \in S$ on $\faktor{Q}{\approx}$. We remark that $T$ is not necessarily a quotient semigroup of $S$, but is in the case that $(Q,S)$ is a transformation semigroup of total functions. 

A congruence $\approx$ is called {\em injective} if every $s \in S$ defines a partial 1-1 function on $\faktor{Q}{\approx}$. It is easy to see that the intersection of injective congruences is injective. Therefore, there is a unique minimal injective congruence $\tau$ on any transformation semigroup $(Q,S)$. We call $\tau$ the Tilson congruence on a transformation semigroup because of the following proved in \cite{Redux}. It is central to the theory of flows.

\bt\label{redux}

Let $(G \times B,S)$ be a $\GM$ transformation semigroup. Then the minimal injective congruence $\tau$ on $S$ is defined as follows: 
$(g,b) \tau (g',b')$ if and only if there are elements $s,t \in S_{II}$ such that $(g,b)s=(g',b')$ and $(g',b')t=(g,b)$.

\et
 
\brm

\item{We can state this theorem by $(g,b)S_{II}=(g',b')S_{II}$. That is, $(g,b)$ and $(g',b')$ define the same ``right coset'' of $S_{II}$ on $G \times B$.}

\item{The proof in \cite{Redux} works on an arbitrary transitive transformation semigroup. We stated it in the case of $\GM$ transformation semigroups because that's how we will use it in this document.}

\item{Theorem \ref{redux} can be used to greatly simplify the proof in \cite{lowerbounds2} for decidability of membership in $S_{II}$ for regular elements of an arbitrary  semigroup.}

\erm

For later use we record the following corollary to Theorem \ref{redux}.

\bc

Let $(G \times B,S)$ be a $\GM$ transformation semigroup. Then $(g,b) \tau (g',b')$ if and only if there are elements 
$s,t \in S_{II} \cap I(S)$ such that $(g,b)s=(g',b')$ and $(g',b')t=(g,b)$. That is, we can choose the elements $s$ and $t$ in Theorem \ref{redux} to be in the 0-minimal ideal of $S$.

\ec

\proof The condition is sufficient by Theorem \ref{redux}. Conversely, assume that $(g,b)\tau (g',b')$. Then there are elements $s,t \in S_{II}$ such that $(g,b)s=(g',b')$ and $(g',b')t=(g,b)$. Since $I(S)$ is a 0-simple semigroup, there are idempotents $e,f \in I(S)$ such that 
$(g,b)e=(g,b)$ and $(g',b')f=(g',b')$. Therefore, $(g,b)es=(g,b)s=(g',b')$ and similarly $(g',b')ft =(g,b)$. Since $e,f$ are idempotents we have $es,ft \in S_{II}$. As $I(S)$ is a 0-minimal ideal, we also have $es,ft \in I(S)$. \qed

%
%

\subsection{Definition of Flows and Their Properties}
Let $(G \times B,S)$ be the transformation semigroup associated to a $\GM$ semigroup $S$ and let $X$ be a generating set for $S$. By a deterministic automaton we mean an automaton such that each letter defines a partial function on the state set.

\bd

Let $\mathcal{A}$ be a deterministic automaton with state set $Q$ and alphabet $X$. 
A {\em flow} to the lattice $SP(G\times B)$ on $\mathcal{A}$ is a function $f:Q \rightarrow SP(G \times B)$ such that for each $q \in Q, x \in X$, with $qf =(Y,\pi)$ and $(qx)f = (Z,\theta)$ we have:

\begin{enumerate}

\item{For all $(g,b) \in G \times B$, there is a $q \in Q$ such that $(g,b) \in qf$.}
 
\item{$Yx \subseteq Z$.}

\item{Multiplication by $x$ considered as an element of $S$ induces a one-to-one partial map from $\faktor{Y}{\pi}$ to $\faktor{Z}{\theta}$.

That is, for all $y,y' \in Y$ we have $y,y'$ are in a $\pi$-class if and only if $yx,y'x$ are in a $\theta$-class whenever $yx,y'x$ are both defined.}

\item{{\bf The Cross Section Condition:} For all $q \in Q$, $qf$ is a cross-section. That is, for all $g,h \in G, b \in B$ if $(g,b),(h,b) \in qf$ it follows that $g=h$.}

\end{enumerate}
\ed


\bd

Let $\mathcal{A}$ be a deterministic automaton with state set $Q$ and alphabet $X$. 
A {\em flow} to the lattice $Rh_{B}(G)$ is a flow to $SP(G\times B)$ such that for each state $q$, $qf \in CS(G\times B)$. That is, $qf$ is a $G$-invariant cross-section.
   
\ed

\brm

    It follows from the original statement of the Presentation Lemma \cite{AHNR.1995} that if $S$ has a flow with respect to some automaton over $SP(G\times B)$, then it has a flow from the same automaton over $Rh_{B}(G)$. The proof of the equivalence of the Presentation Lemma and Flows in Section 3 of \cite{Trans} works as well for the Presentation Lemma in the sense of \cite{AHNR.1995}. 

\erm

We will use all the terminology and concepts from \cite{Trans}. If $\mathcal{A}$ is an automaton with alphabet $X$ and state set $Q$, recall that its completion is the automaton $\mathcal{A}^{\square}$ that adds a sink state $\square$ to $Q$ and declares that $qx = \square$ if $qx$ is not defined in $\mathcal{A}$. A flow on $\mathcal{A}$ is a complete flow if on $\mathcal{A}^{\square}$ extending $f$ by letting $\square f = (\emptyset, \emptyset)$ the bottom of both the set-partition and Rhodes lattices remains a flow and furthermore, for all $(g,b) \in G \times B$, there is a $q\in Q$ such that $(g,b) \in Y$, where $qf=(Y,\Pi)$. All flows in this paper will be complete flows.

 Flows are related to the Presentation Lemma \cite{AHNR.1995}, \cite[Section 4.14]{qtheory}. They give a necessary and sufficient condition for $Sc =\operatorname{RLM}(S)c$ where $S$ is a $\operatorname{GM}$ semigroup.  See Section 3 of \cite{Trans} for a proof of the following Theorem

\bt\label{PLflow}[The Presentation Lemma-Flow Version]

Let $(G \times B,S)$ be a $GM$ transformation semigroup with $S$ generated by $X$. Let $k> 0$ and assume that $\operatorname{RLM}(S)c = k$.  Then $Sc = k$ if and only if there is an $X$ automaton $\mathcal{A}$ whose transformation semigroup $T$ has complexity strictly less than $k$ and a complete flow $f:Q \rightarrow Rh_{B}(G)$.


\et

 Flows were preceded by the Presentation Lemma \cite{AHNR.1995}, \cite[Section 4.14]{qtheory}. The Presentation Lemma was shown to be equivalent to the existence of an appropriate flow in \cite[Section 3]{Trans}. There are three versions of the Presentation Lemma and its relation to Flow Theory in the literature \cite{AHNR.1995}, \cite[Section 4.14]{qtheory}, \cite{Trans}. These have very different formalizations and it is not clear how to pass from one version to another. The terminology is different as well. For example, the definition of cross-section in each of these references, as well as the one we use in this paper is different.  The following Theorem gives a unified approach to these topics. It summarizes known results in the literature and is meant to emphasize the strong connections between slices, flows and the presentation lemma and the corresponding direct product decomposition. \cite{FlowsI} for a proof. We use some of the results that we've proved in previous sections. For background on the derived transformation semigroup and the derived semigroup theorem see \cite{Eilenberg, qtheory}.

\bt \label{uniform}

Let $(G \times B,S)$ be $\GM$ and assume that $\RLM(S)c \leq n$. Then the following are equivalent:

\begin{enumerate}

\item{$Sc \leq n$.}

\item{There is an aperiodic relational morphism $\Theta:S \rightarrow H\wr T$, where $H$ is a group and $Tc \leq n-1$.}

\item{There is a relational morphism $\Phi:S \rightarrow T$ where $Tc \leq n-1$ and such that the Derived Transformation semigroup $D(\Phi)$ is in $Ap*Gp$.}

\item{There is a relational morphism $\Phi:S \rightarrow T$ where $Tc \leq n-1$ such that the Tilson congruence $\tau$ on the Derived Transformation semigroup $D(\Phi)$ is a cross-section.}

\item{$(G \times B,S)$ admits a flow from a transformation semigroup $(Q,T)$ with $(Q,T)c \leq n-1$.}

\item{$S \prec (G \wr (Sym(B)) \wr T) \times \RLM(S)$ for some transformation semigroup $T$ with $Tc \leq n-1$.}

\end{enumerate}

\et

%
%

We emphasize the case for complexity 1 and aperiodic flows as this is used in many of our examples.

\bt \label{Apuniform}

Let $(G \times B,S)$ be $\GM$ and assume that $\RLM(S)c \leq 1$. Then the following are equivalent:

\begin{enumerate}

\item{$Sc=1$.}

\item{There is an aperiodic relational morphism $\Theta:S \rightarrow H\wr T$, where $H$ is a group and $T$ is aperiodic.}

\item{There is a relational morphism $\Phi:S \rightarrow T$ where $T$ is aperiodic such that the Derived transformation semigroup $D(\Phi)$ is in $Ap*Gp$.}

\item{There is a relational morphism $\Phi:S \rightarrow T$ where $T$ is aperiodic such that the Tilson congruence $\tau$ on the Derived transformation semigroup $D(\Phi)$ is a cross-section.}

\item{$(G \times B,S)$ admits an aperiodic flow.}

\item{$S \prec (G \wr (Sym(B)) \wr T) \times \RLM(S)$ for some aperiodic semigroup $T$.}

\end{enumerate}

\et

\section{APPENDIX The Flow Monoid and the Evaluation Transformation Semigroup}\label{Eval}

\subsection{The Monoid of Closure Operations}\label{Closops}

In this Appendix we gather definitions and results from \cite{Trans}. For more details the reader is strongly urged to consult this paper. Since all the computations we do on examples in this document are done in the Evaluation Transformation Semigroup defined in Section 5 of \cite{Trans} our goal is to summarize the background material in that paper needed to define this object. We begin with the definition of the monoid of closure operations on the direct product $L^{2} = L \times L$ of a lattice $L$ with itself. 

Let $L$ be a lattice and let $L^{2} = L \times L$. Let $f$ be a closure operator on $L^2$. By definition this means that $f$ is an order preserving, extensive (that is, for all $(l_{1},l_{2}) \in L^{2}, (l_{1},l_{2}) \leqslant (l_{1},l_{2})f$), idempotent function on $L^2$. A {\em stable pair} for $f$ is a closed element of $f$. Thus a stable pair $(l,l') \in L^2$ is an element such that $(l,l')f = (l,l')$.  The stable pairs of $f$ are a meet closed subset of $L \times L$. Conversely, each meet closed subset of $L \times L$ is the set of stable pairs for a unique closure operator on $L \times L$. We will identify $f$ as a binary relation on $L$ whose pairs are precisely the stable pairs. Let $B(L)$ be the monoid of binary relations on the set $L$.  As is well known, $B(L)$ is isomorphic to the monoid $ M_{n}(\mathcal{B})$ of $n \times n$ matrices over the 2-element Boolean algebra $\mathcal{B}$, where $n=|L|$. The Boolean matrix associated to $f$ is of dimension $|L| \times |L|$. It has a 1 in position $(l_{1},l_{2})$ if $(l_{1},l_{2})$ is a stable pair and a 0 otherwise. 

With this identification, the collection $\mathcal{C}(L^{2})$ of all closure operators on $L^2$ is a submonoid of the monoid $B(L)$ of binary relations on $L$~\cite[Proposition 2.5]{Trans}. We thus also consider $\mathcal{C}(L^{2})$ to be a monoid of $|L| \times |L|$ Boolean matrices. Let $L$ be either $Rh_{B}(G)$ or $\operatorname{SP}(G\times B)$. We now recall the definition of some important unary operations on $\mathcal{C}(L^{2})$. 
\bd

Let $f \in \mathcal{C}(L^{2})$.

\begin{enumerate}

\item{The domain of $f$ denoted by $\operatorname{Dom}(f)$ is the set $\{x \mid \exists y, (x,y) \in f\}.$}

\item{Define the relation $\overleftarrow{f}$ by $\overleftarrow{f}= \{(x,x) \mid x \in \operatorname{Dom}(f)\}$. $\overleftarrow{f}$ is called {\em back flow along $f$}. See \cite[Remark 2.26]{Trans} for the reason for this terminology.}

\item{The relation $f^{*}$ is defined by $f^{*} = f \cap \{(x,x) \mid x \in X\}$. $f^{*}$ is called {\em the Kleene closure of $R$}. See Section 2 and Section 4 of \cite{Trans} for the reason for this terminology.}

\item{Define the {\em loop of $f$} to be the relation $f^{\omega+*}=f^{\omega}f^{*}$, where $f^{\omega}$ is the unique idempotent in the subsemigroup generated by $f$.}

\end{enumerate}

\ed

At times it is convenient to identify $\overleftarrow{f}$ as $1|_{\operatorname{Dom}(f)}:L \rightarrow L$, the identity function restricted to 
$\operatorname{Dom}(f)$. Similarly, we identify $f^{*}$ as $1|_{\operatorname{Fix{f}}}:L \rightarrow L$, the restriction of the identity to the set of fixed-points of $f$, where $x \in L$ is a fixed point if $(x,x) \in f$. Our use of these will be clear from the context.

\subsection{The 0-Flow Monoid}\label{FMon}

Let $S$ be a $\GM$ semigroup generated by $X$ and let $x \in X$. Let $I(S) = \mathcal{M}^{0}(A,G,B,C)$. We define a binary relation $f_{x}$
on $\operatorname{SP}(G \times B)$ by $((Y,\Pi), (Z,\Theta)) \in f_{x}$ if and only if $Yx \subseteq Z$ and the partial function induced by
right multiplication by $x$, $\cdot x: Y \rightarrow Z$ induces a well-defined partial injective map $\cdot x: \faktor{Y}{\pi} \rightarrow Z/\Theta$. This means that if  $(g,b), (g',b') \in Y$ and $(g,b)x,(g',b')x$ are both defined (and hence in $Z$ ), then $(g,b)\Pi (g',b')$ if and only if $(g,b)x\Theta (g',b')x$. Then $f_{x} \in \mathcal{C}(L^{2})$. See \cite[Proposition 2.22]{Trans}. $f_{x}$ is called the free-flow by $x$.

We now define the 0-flow monoid $M_{0}(L)$ as follows.

\bd\label{Flowops}

Let $L = \operatorname{SP}(G \times B)$. The 0-flow monoid $M_{0}(L)$, is the
smallest subset of $\mathcal{C}(L^{2})$ satisfying the following axioms:

\begin{enumerate}

\item {(Identity) The multiplicative identity $I$ of $\mathcal{C}(L^{2})$ is in $M_{0}(L)$.}

\item{(Points) For all $x \in X$, $f_{x}$ the free-flow along $x$ belongs to $M_{0}(L)$.}

\item{(Products) If $f_{1}, f_{2} \in M_{0}(L)$, then $f_{1}f_{2} \in M_{0}(L)$.}

\item{(Vacuum) If $f \in M_{0}(L)$, then $\overleftarrow{f} \in M_{0}(L)$.}\label{Vac}

\item{(Loops) If $f \in M_{0}(L)$, then $f^{\omega+*}\in M_{0}(L)$.}

\end{enumerate}

\ed

We remark that for each $n \geq 0$, there is an $n$-flow monoid $M_{n}(L)$ defined in \cite{Trans}. The definition above of $M_{0}(L)$ is exactly what is defined in \cite{Trans} and used extensively in \cite{complexity1}. For $n>0$, Axioms (1)-(4) are the same for $M_{n}(L)$ as for $M_{0}(L)$. Axiom (5) restricts the use of the loop operator to $n$-loopable elements \cite[Section 4]{Trans}. In \cite{complexityn}, $n$-loopable elements are replaced by a more restrictive definition and the modified $M_{n}(L)$ plays a crucial role in the main results of \cite{complexityn}.

\subsection{The Evaluation Transformation Semigroup}\label{Ets}

We now defined the Evaluation Transformation Semigroup $\mathcal{E}(L) = (\operatorname{States}, \operatorname{Eval}(L))$ as in \cite{Trans}. Again, because of our interest in this paper on aperiodic flows, we don't define the $n$-Evaluation Transformation Semigroups $E_{n}$ for $n>0$. We begin with the definition of Well-Formed Formualae (WFFs).

\bd

Let $X$ be an alphabet. We define a well-formed formula inductively as follows.

\begin{enumerate}

\item{The empty string $\epsilon$ is a well-formed formula.}

\item{Each letter $x \in X$ is a well-formed formula.}

\item{If $\tau, \sigma$ are well-formed formulae, then
so is $\tau\sigma$.}

\item{If $\tau$ is a well-formed formula that is not a proper power (i.e., not of the form
$\sigma^{n}$ where $n > 1$), then $\tau^{\omega+*}$ is also a well-formed formula.} 

\end{enumerate}

\ed

The set of well-formed
formulae is denoted by $\Omega(X)$. Well-formed formulae will be denoted by Greek
letters. As a convention, if $\tau=\sigma
^{n}$, where $\sigma$ is not a proper power, then we set
$\tau^{\omega+*}=\sigma^{\omega+*}$. In other words, we extract roots before applying the unary operation $\omega+*$.

Let $V = \prod_{f \in M_{0}(L)}\overleftarrow{f}$. $V$ is called the {\em Vacuum}. See \cite{Trans} for an explanation of this terminology. In \cite{Trans}, $V$ is denoted by $\mathcal{F}_{0}$. It is proved in \cite{Trans} that $V$ is an idempotent in $M_{0}(L)$. We now will work exclusively in the subsemigroup $VM_{0}(L)V$ of $M_{0}(L)$. We want to interpret WFFs in $VM_{0}(L)V$.

\bd

Define recursively a partial function $\mathcal{I}: \Omega(X)\rightarrow VM_{0}(L)V$ as follows. 

\begin{enumerate}

\item{$\epsilon\mathcal{I} = V$.}

\item{$x\mathcal{I} = VxV$ for $x \in X$.}

\item{If $\mathcal{I}$ is already defined on $\tau, \sigma \in \Omega(X)$, set $(\tau\sigma)\mathcal{I} = \tau\mathcal{I}\sigma\mathcal{I}$.}

\item{If $\tau \in \Omega(X)$ is not a proper power and $\tau\mathcal{I}$ is defined, set $\tau^{\omega+*}\mathcal{I} = 
(\tau\mathcal{I})^{\omega+*}$.}

\end{enumerate}

\ed

We normally omit $\mathcal{I}$ and assume that a WFF $\tau$  is being evaluated in $VM_{0}(L)V$ according to the definition of $\mathcal{I}$. 
 We first define a new operator on elements of $M_{0}(L)$ called {\em forward flow}. Recall that the bottom of the lattice $\operatorname{SP}(G \times B)$ is the pair $\Box = (\emptyset, \emptyset)$.

\bd

Let $f \in M_{0}(L)$ and let $l \in L$. Let $(l,\Box)f =(l_{1},l_{2})$ We define the forward flow of $f$ denoted by $\overrightarrow{f}$ by 
$l\overrightarrow{f}=l_{2}$. That is, we apply $f$ to $(l,\Box)$ and project to the right-hand coordinate.

\ed

The following is proved in Sections 2 and 4 of \cite{Trans}.

\begin{enumerate}

\item{$\overrightarrow{f}:L \rightarrow L$ is an order preserving function on $L$.}

\item{If $lV=l$, then $(l,\Box)f=(l,l')$ and $l' \in LV$. That is, the left-coordinate of $(l,\Box)f$ is still $l$ and the right-hand coordinate is in $LV$. Therefore $\overrightarrow{f}$ is a well-defined function from $LV$ to $LV$.}

\item{The assignment of $f$ to $\overrightarrow{f}$ defines an action of $VM_{0}(L)V$ on $LV$. It follows that we have a transformation semigroup $(LV,M'_{0}(L))$, where $M'_{0}(L)$ is the image of $VM_{0}(L)V$ on $LV$ under this action.}

\end{enumerate}

We now restrict the action in $(LV,M'_{0}(L))$ to the set of $\operatorname{States}$ defined as follows. Let $(g,b) \in G \times B$. Then 
the element $(\{(g,b)\},\{(g,b)\}) \in L$ is called a {\em point}. In Section 5 of \cite{Trans}, it is proved that every point $p$ satisfies $pV=p$. 

\bd

The set of $\operatorname{States}$ is the smallest subset of $LV$ such that:

\begin{enumerate}

\item{Every point $p \in \operatorname{States}$.}

\item{If $l \in \operatorname{States}$, then $l\overrightarrow{f} \in \operatorname{States}$.}

\end{enumerate}

\ed

In other words $\operatorname{States}$ is the smallest subset of $LV$ containing the points and closed under the action of $M_{0}(L)$ on $LV$. We can finally define the Evaluation Transformation Semigroup where all the computations in the examples in this paper take place.

\bd

The Evaluation Transformation Semigroup $\mathcal{E}(L)$ is defined by $\mathcal{E}(L) = (\operatorname{States}, \operatorname{Eval}(L))$ where $\operatorname{Eval}(L)$ is the image of $M'_{0}(L)$ by restricting its action to $\operatorname{States}$.

\ed

We remark that there is an Evaluation Transformation Semigroup $\mathcal{E}_{n}(L)$ for all $n \geq 0$. The definition here is the case $n=0$
which is all we need in this paper as we are only concerned with aperiodic flows.

\bibliography{stubib}
\bibliographystyle{abbrv}

{(Stuart Margolis) Department of Mathematics, Bar Ilan University, Ramat Gan 52900, Israel

{\it Email address}: \; \texttt{margolis@math.biu.ac.il}}

\medskip

{(John Rhodes) Department of Mathematics, University of California, Berkeley, CA 94720, U.S.A.

{\it Email address}: \; \texttt{blvdbastille@gmail.com}}

\end{document}